\documentclass[11pt,twoside,reqno,draft]{article}
\usepackage[centertags]{amsmath}
\usepackage{amssymb,amsfonts,amsthm}
\newcommand{\TryPackage}[3]{\IfFileExists{#1.sty}{\usepackage{#1} #2}{#3}}
\TryPackage{mathrsfs}{\renewcommand{\mathcal}{\mathscr}}{%
   \TryPackage{eucal}{}{}}

\usepackage{a4}             
\usepackage{math40,local}   

\newcommand{\signature}{\relax}

\begin{document}
\title{On the $L^2$--Stokes theorem and Hodge theory for
singular algebraic varieties}
\author{Daniel Grieser and Matthias Lesch}
\date{\today}
\renewcommand{\signature}{\bigskip\noindent
Humboldt Universit\"at zu Berlin\\
Institut f\"ur Mathematik\\ 
Sitz: Rudower Chaussee 25\\
D--10099 Berlin\\[1em]
Universit\"at zu K\"oln,\\
Mathematisches Institut,\\
Weyertal 86-90,\\
D--50931 K\"oln
\bigskip\noindent
\begin{tabbing}
Email--addresses: \=
   grieser@mathematik.hu-berlin.de\\
\>lesch@math.uni-koeln.de\\[1em]
Internet: \> http://mi.uni-koeln.de/$\sim$lesch/\\
\> http://www.mathematik.hu-berlin.de/Math-Net/members/grieserd.rdf.html
\\
\end{tabbing}
}
\maketitle
\begin{abstract}
For a projective
 algebraic variety $V$ with isolated singularities, endowed with a
metric induced from an embedding, we consider the analysis of the natural
partial differential operators on the regular part of
 $V$. We show that, in the complex case,
the Laplacians of the de Rham and Dolbeault complexes
are discrete operators except possibly in degrees $n,n\pm 1$, where $n$ is the
complex dimension of $V$.
We also prove a Hodge theorem on the operator level
and the $L^2$--Stokes theorem outside the degrees $n-1,n$.
We show that the $L^2$-Stokes theorem may fail to hold
in the case of real algebraic varieties, and also discuss the
$L^2$-Stokes theorem on more general non-compact spaces.

\end{abstract}
\tableofcontents
%


\section{Introduction}
The interplay between geometric 
differential operators on a Riemannian manifold
 and the geometry of the underlying manifold has been the focus of many
 efforts; one of 
the early highlights is the Atiyah-Singer Index Theorem.

Since the work of Atiyah and Singer one has become more and more interested in 
various types of manifolds with singularities. While the case of a smooth
compact
manifold is fairly well understood, 
the picture is far from  complete for singular manifolds.
It is impossible to give a 
complete account of the existing literature here. We only mention Cheeger's 
work on manifolds with cone-like singularities \cite{Che:SGSRS}, 
Melrose's b-calculus \cite{Mel:APSIT} and Schulze's calculus on
singular manifolds \cite{Schul:PDBVPCSA}.

However, singularities occuring in "the real world" are often much more
complicated than just conical. 
A very natural class of singular spaces is the class of 
(real or complex) projective algebraic varieties. These 
are special cases of stratified spaces. Topologically, 
stratified spaces are of iterated cone-type, and probably this was 
Cheeger's main motivation to develop an analysis of elliptic 
operators on such manifolds. However, it seems that the inductive step, 
i.e. the generalization 
of Cheeger's theory to stratified spaces,  has still not been done. 
A more serious problem is that the natural metrics on algebraic varieties, 
i.e. those induced from a metric on projective space are not 
of iterated cone--type. A great deal of efforts have been made to find local 
models of such metrics \cite{Gri:LGSRAS}, \cite{HsiPat:LCNAS}, \cite{Pat:LATIS}.

Nevertheless, there exist partial results about the interplay 
between $L^2$--cohomology and intersection cohomology \cite{Ohs:LCCS},
\cite{Ohs:LCCSII}, 
mixed Hodge structure \cite{ParSte:PHSLCVIS},
and the so-called $L^2$-K\"ahler package \cite{Che:HTCC}, \cite{BruLes:KHTCCC},
\cite{ParSte:LCCPV}.
Having the $L^2$--K\"ahler package on a complex 
algebraic variety in general would be extremely nice, because it implies many 
of the fundamental operator identities which one has in the 
compact case. 

The problems we discuss in this note are the $L^2$--Stokes theorem and 
discreteness of the Laplace-Beltrami operators:

Let $V$ be a real or complex projective variety and let 
$M:=V\setminus \sing V$ be its regular part. We equip $M$ with the 
Riemannian metric $g$ induced by a smooth metric on projective space (in the complex 
case we assume $g$ to be K\"ahler). Furthermore, let $(\Omega_0(M),d)$ 
be the de Rham complex of differential forms acting on smooth forms 
with compact support. A priori, the operator $d$ has several closed 
extensions in the Hilbert space of square integrable forms. These lie
between the "minimal" and the "maximal" one. The latter are defined by 
\begin{align*}
d_{\min} & := \overline{d}  = \text{closure of } d, \\
d_{\max} & := (d^t)^*  = \text{adjoint of the formal adjoint } d^t \text{ of } d.
\end{align*}
Some authors address the maximal as the Neumann and 
the minimal as the Dirichlet extension. We do not adopt this terminology 
since it may be misleading: On a compact manifold with boundary the 
Laplacians corresponding to the  maximal/minimal extensions both are of 
mixed Dirichlet/Neumann type. The maximal/minimal extensions of $d$ 
produce so-called Hilbert complexes. A detailed account of the 
functional analysis of Hilbert complexes was given in \cite{BruLes:HC}. 
We note that the cohomology of the $d_{\max}$ complex is the celebrated 
$L^2$--cohomology 
$$ H^i_{(2)}(M):=\ker d_{i,\max}/\im d_{i-1,\max},$$
which has been the subject of intensive studies.

Having to make the distinction between $d_{\max}$ and $d_{\min}$ can be 
tedious.
If $d_{\max}\not=d_{\min}$ then even simple facts known from compact manifolds
might not be true. Therefore, it is desirable
to have \emph{uniqueness}, i.e. $d_{\max}=d_{\min}$, which is better known as
the \emph{$L^2$--Stokes theorem (\LST)} because of its equivalent formulation
(\ref{eq6}). We would like to emphasize that validity of \LST\ does \emph{not} imply
essential selfadjointness of the Laplace-Beltrami operator $\Delta$
(defined on compactly
supported smooth forms), as can be seen in the case of cones already. Instead,
 \LST\ is equivalent to the selfadjointness of the specific extension
$d^t_{\min}d_{\min} + d_{\min}d^t_{\min}$ of $\Delta$.

On a compact manifold without boundary the 
\LST\ follows by a simple mollifier argument. 
It is known to hold 
for several types of non-compact manifolds. We give a short 
account of known results and their proofs
 in Section 2. The conjecture that the \LST\
holds for complex projective varieties has been around implicitly for 
quite a while, e.g. in \cite{CheGorMac:LCIHSAV}, although the authors cannot give a 
reference where it is stated explicitly, except for complex surfaces 
\cite{Nag:RLCSAS}. 

The \LST\ is quite plausible for complex projective varieties
because all strata are of even codimension, in particular there is no
boundary. So, no "boundary terms" should occur in the integration by parts
implicit in the \LST. However, this picture cannot be 
complete since we will prove the following (Proposition \plref{2.10}).

\begin{theorem}\label{S-intro.1} 
There exists an even dimensional real projective 
variety whose singular set consists of a single point, 
such that the $L^2$--Stokes 
theorem does not hold. 
\end{theorem}

Thus, if $L^2$ST holds for complex projective varieties the reason must lie 
in the complex structure. 

Call two Riemannian metrics $g$, $\tilde{g}$ on $M$ \emph{quasi-isometric}
if there is a constant $C$ such that for each $x\in M, v\in T_xM$ one has
$$ C^{-1} g_x(v,v) \leq \tilde{g}_x(v,v) \leq C g_x(v,v).$$
The domains of $d_{\min}$ and $d_{\max}$ and therefore 
validity of \LST\ are quasi-isometry invariants. Therefore, if \LST\ 
holds for one metric
induced from projective space, it holds for all such metrics. Thus its 
validity is independent of the K\"ahler structure. But surprisingly enough 
the K\"ahler structure can be very useful as a tool. For example the 
estimate \cite[Prop. 2.13]{ParSte:PHSLCVIS} (see also Proposition 
\ref{3.1} below)
proved by Pardon and Stern on complex projective varieties with 
isolated singularities makes essential use of the K\"ahler structure. 
They apply their estimate to derive a Hodge structure on the $L^2$--cohomology. 
In particular, 
they conclude that in various degrees the cohomologies of the $d_{\min}$ and 
$d_{\max}$ de Rham complexes coincide. Together with results in
\cite{BruLes:KHTCCC} the estimates of Pardon and Stern can be used to prove more:
\begin{theorem}\label{3.3}
Let $V\subset \CP^N$ be an algebraic variety with isolated 
singularities, of complex dimension $n$, and let
$M=V\setminus \sing V$, equipped with a K\"ahler 
metric induced by a K\"ahler metric on $\CP^N$.
Then, for $k\not= n-1,n$ resp.\ $p+q\not= n-1,n$, the $L^2$ Stokes theorem
holds and we have uniqueness for the Dolbeault operators, i.e.
\begin{align}
d_{k,\max} & = d_{k,\min}, \quad\quad k\not= n-1,n, \label{3a}\\
\pl_{p,q,\max} & = \pl_{p,q,\min}, \quad\quad p+q \not= n-1,n. \label{3b}
\end{align}

Furthermore, for $k\not= n,n\pm 1$ we have
\begin{equation}
\begin{split}\label{4}
d_{k-1,\min} &d^t_{k-1,\min} + d^t_{k,\min} d_{k,\min}\\
=& \Delta^{\cF}_k  = 2 \bigoplus_{p+q=k} \Delta^{\cF}_{p,q,\pl} \\
=& 2 \bigoplus_{p+q=k} \pl_{p-1,q,\min} \pl^t_{p-1,q,\min} +
                       \pl^t_{p,q,\min} \pl_{p,q,\min}, 
\end{split}                       
\end{equation}
i.e. the Hodge decomposition holds in the operator sense.
\end{theorem}
Here, $\Delta_k^{\cF}$ is the Friedrichs extension, see Section 3.
We expect that $L^2$ST is true in fact for all degrees and that
the first equality in (\ref{4}) holds in all degrees except
$n$ (in degree $n$ it will usually not hold).
This would be a very interesting result since it would imply
the K\"ahler package, as shown in \cite[Th. 5.8]{BruLes:KHTCCC}.

As already mentioned, essential self-adjointness cannot be expected for the 
Laplacian. It is therefore quite surprising to obtain its Friedrichs 
extension from the de Rham complex. It was shown in 
\cite{BruLes:KHTCCC} that this case is exeptional and that 
it has some nice consequences.

It is interesting to know more about the 
structure of the spectrum of the Friedrichs extension of the Laplacian. 
We will prove the following.
\begin{theorem}\label{3.5} Under the assumptions of Theorem \plref{3.3}
the operator on $k$-forms
$$ \Delta^{\cF}_k=d_{k-1,\min} d^t_{k-1,\min} + d^t_{k,\min} d_{k,\min}$$
is discrete for $k\not=n,n\pm1$.
\end{theorem}
For $k=0$ on an algebraic surface or threefold with
isolated singularities, this follows from \cite{Nag:HONSAS} and \cite{Pat:HTSAT}, 
where an estimate for the
heat kernel is proved. For algebraic curves, the full asymptotic expansion
of the heat trace was proved by Br\"uning and Lesch \cite{BruLes:SGAC}.
The existence of such an asymptotic expansion in
general remains a challenging open problem.

This note is organized as follows: In Section 2 we discuss various 
aspects of the $L^2$ Stokes theorem and its history. Furthermore, we prove 
Theorem \plref{S-intro.1}. Section 3 is devoted to complex projective 
varieties and the proof of Theorems \ref{3.3} and \ref{3.5}.

\bigskip
\centerline{\bf Acknowledgement}

Both authors were supported by the Gerhard-Hess program of 
Deutsche Forschungsgemeinschaft.

\section{The $L^2$--Stokes theorem}

We start with some general remarks about elliptic complexes and 
their ideal boundary conditions on non-compact manifolds. In 
particular we recall some results from \cite{BruLes:HC,BruLes:KHTCCC}.

Let $(M,g)$ be a Riemannian manifold and 
\begin{equation*} 
D:C_0^\infty(E) \longrightarrow C_0^\infty(F)
\end{equation*}
a first order differential operator between sections of the hermitian vector
bundles $E,F$. We consider $D$ as an unbounded operator $L^2(E)\to L^2(F)$
and define two closed extensions of $D$ by
\begin{align*}
D_{\min} & := \overline{D}  = \text{closure of } D, \\
D_{\max} & := (D^t)^*  = \text{adjoint of the formal adjoint } D^t \text{ of } D.
\end{align*}
Note that $D_{\min} \subset D_{\max}$ and $D_{\min} = (D^*)^* = (D^t_{\max})^*$
where we write $D^t_{\max} = (D^t)_{\max}$. The domains of
$D_{\maxmin}$ can be described as follows:
\begin{equation}
\begin{split} \label{eq3}
\dom(D_{\min}) & =  \Big\{s\in L^2(E)\,\Big|\, 
     \begin{array}{l} \text{ There exists a sequence } (s_n)\subset 
                C_0^\infty(E) \\
       \text{ with } s_n\to s, Ds_n \to Ds \text{ in } L^2 (E)
     \end{array} \Big\},\\
\dom(D_{\max}) & =  \{s\in L^2(E)\,|\, Ds\in L^2(E)\}.
\end{split}
\end{equation}
$D_{\max}$ is maximal in the sense that it does not have a proper closed
extension that has $C_0^\infty(F)$ in the domain of its adjoint.
Of course, $D_{\max}$ does have (abstractly defined) proper closed
extensions without this property.

The following well-known fact shows that the 'difference' between $D_{\max}$ and $D_{\min}$
only depends on the behavior of $(M,g)$ at 'infinity', i.e.\ when leaving
any compact subset of $M$. We include a proof for completeness.

\begin{lemma}\label{lemma2.1}
If $s\in\dom(D_{\max})$ and $\phi\in C_0^\infty(M)$ then 
$\phi s\in \dom(D_{\min})$.
\end{lemma}

\begin{proof}
We use a Friedrichs mollifier, i.e. a family of operators
$J_\eps:\ce'(M) \to C_0^\infty(M), \eps\in (0,1],$ such that
$J_\eps f\to f$ in $L^2(M)$ for any $f\in L^2_{\text{comp}}(M)$ and such 
that $J_\eps$ and the commutator $[D,J_\eps]$ are bounded operators
on $L^2(M)$, uniformly in $\eps$. For the existence of such mollifiers
see \cite[Ch.\ II.7]{Tay:PO}.

Now we have $J_\eps(\phi s)\in \dom(D_{\min})$,
$J_\eps(\phi s)\to\phi s$, and $D(J_\eps(\phi s)) = 
[D,J_\eps](\phi s) + J_\eps (D(\phi s))$ is uniformly bounded in $L^2(M)$,
as $\eps\to 0$.
Therefore, there is a constant $C$ such that for all $t\in\dom(D^*)$ we
have
\begin{equation}\label{domest}
|(\phi s, D^*t)| = \lim_{\eps\to0} |(J_\eps(\phi s),D^*t)|
 = \lim_{\eps\to0} |(D(J_\eps(\phi s)),t)| \leq C\|t\|. 
\end{equation}
This means $\phi s\in \dom((D^*)^*) = \dom(D_{\min})$.
\end{proof}

We now turn to elliptic complexes. Elliptic complexes on manifolds with
singularities have been studied systematically e.g.\ in
\cite{RemSchul:ITEBP},
\cite{Schul:ECMWCS}. For a general discussion of Fredholm complexes and
Hilbert complexes we also refer to \cite{Seg:FC} and \cite{BruLes:HC}.
 Let
\begin{equation}\label{eq4}
 (C_0^\infty(E),d): 0\longrightarrow \cinfz{E_0} \stackrel{d_0}{\longrightarrow} \cinfz{E_1}
\stackrel{d_1}{\longrightarrow} \ldots \stackrel{d_{N-1}}{\longrightarrow} \cinfz{E_N} \longrightarrow 0 
\end{equation}
be an elliptic complex. The main examples are the de Rham complex
$(\Omega^*_0(M),d)$ and, for a K\"ahler manifold, the Dolbeault complexes
$(\Omega^{*,q}_0 (M), \pl_{\cdot,q})$ and $(\Omega^{p,*}_0(M), \plbar_{p,\cdot})$.

Next we recall the notion of {\em Hilbert complex} (cf.\ \cite{BruLes:HC},
for example).
This is a complex 
\begin{equation}\label{eq5}
(\dom,D): 0 \longrightarrow \dom_0 \stackrel{D_0}{\longrightarrow} \dom_1 \stackrel{D_1}{\longrightarrow}
\ldots \stackrel{D_{N-1}}{\longrightarrow} \dom_N \longrightarrow 0
\end{equation}
of closed operators $D_k$ with domains $\dom_k$ lying in Hilbert spaces $H_k$.
An {\em ideal boundary condition (ibc)} of an 
elliptic complex $(\cinfz{E},d)$ is a choice of extensions $D_k$
of $d_k$ which form a Hilbert complex (with $H_k=L^2(E_k)$).
That is, an ibc is a choice of closed extensions $D_k$
of $d_k, 0\le k<N$, with the additional property that $D_k(\cd_k)\subset
\cd_{k+1}$. $d_{k,\min}$ and $d_{k,\max}$ are examples of ibc's for any
elliptic complex.

The main question that we address here is whether there is only one ibc,
i.e.\ whether $d_{k,\min} = d_{k,\max}$ for all $k$. This is called the 
case of {\em uniqueness}. Uniqueness is equivalent to
\begin{equation}\label{eq6}
 (d_{\max} \omega,\eta) = (\omega,d^t_{\max} \eta) \quad \text{ for all } \omega\in\dom(d_{\max}),
     \eta\in\dom(d^t_{\max}).
\end{equation}     
\begin{dfn}\label{a}
We say that the {\em $L^2$--Stokes theorem} (\LST)
holds for $(M,g)$ if \myref{eq6}\ is true
for the de Rham complex on $M$.
\end{dfn}
This means that no boundary terms appear in the integration by parts that
is implicit in \myref{eq6}. In particular, \LST\ holds for closed $M$.

We now return to general elliptic complexes \myref{eq4}\ and their
ibc's \myref{eq5}. Given an ibc $(\dom,D)$, define the associated Laplacian
by
$$ \Delta (\dom,D) = (D+D^*)^2.$$
This is a self-adjoint operator. The main instances are $\Delta^{a/r}$,
which are associated with $D=d_{\maxmin}$, respectively. Note that $\Delta^a$
and $\Delta^r$ are not comparable unless they are equal. Lemma 3.1 from
\cite{BruLes:KHTCCC}\ says (as a special case) that this happens if and
only if $d_{\max}=d_{\min}$. Explicitly, the restriction of $\Delta^{a/r}$ to
sections of $E_k$ (e.g. $k$-forms) is given by
$$ \Delta^{a/r}_k = d_{k-1,\maxmin} d^t_{k-1,\minmax}
                      + d^t_{k,\minmax} d_{k,\maxmin}. $$
Also, we denote
$$ \Delta = (d+d^t)^2 \quad \text{ on } \cinfz{E}$$
and $\Delta_k$ its restriction to sections of $E_k$.

\begin{prop}\label{prop2.3}
Let $(\cinfz{E},d)$ be an elliptic complex.
\begin{enumerate}
\item $\cinf{E} \cap \dom(d_{\maxmin})\cap\cd(d^t_{\minmax})$ is (graph) dense in 
$\dom(d_{\maxmin})$.
\item Fix $k$ and set
$$\tilde{\Delta}_k := d_{k-1,\min} d^t_{k-1,\min} + d^t_{k,\min} d_{k,\min}.$$
$\tilde{\Delta}_k$ is self-adjoint if and only if
 $\dmaxdmin{k}$ and $\dmaxdmin{k-1}$.

In particular, if $\Delta_k$ is essentially
self-adjoint for all $k$ then uniqueness holds.
\end{enumerate}
\end{prop}
\begin{proof}
1. By ellipticity, we have
\begin{equation} \label{eq8}
\bigcap_{n\geq 0} \dom((\Delta_k^{a/r})^n) \subset C^\infty(E_k).
\end{equation} 
By Lemma 2.11 in \cite{BruLes:HC}\ the left hand side in \myref{eq8}\
is a core for $d_{k,\maxmin}$ (i.e.\ dense in the graph topology), 
so the claim follows.

2. 'If' is obvious. To prove the converse, observe first
$\tilde{\Delta}_k \subset \Delta_k^{a/r}$. Since
all of these operators are self-adjoint, this implies
$\Delta_k^r = \tilde{\Delta}_k = \Delta_k^a$. Since $\dom(\Delta^{a/r}_k)$
is a core for $d_{k,\maxmin}$ and for $d^t_{k-1,\minmax}$ we are done.
\end{proof}

Now we turn to the de Rham complex and the \LST.

First we note that the validity of \LST\ is a quasi-isometry invariant
since the domains $\dom(d_{k,\maxmin})$ are quasi-isometry invariants
in view of their characterization \myref{eq3}. (Note that $d_k$ itself
is independent of the metric.)

We will be mainly interested in the case of projective varieties with
induced metrics. That is, we start with a variety $V$ in real or 
complex projective space and let $M=V\setminus \sing V$ be its regular part.
The metric on $M$ is obtained by restriction of a Riemannian metric
on projective space. Since any two Riemannian metrics on projective
space are quasi-isometric, all metrics on $M$ obtained in this way are mutually
quasi-isometric, so if \LST\ holds for one such metric then it holds for
all.

The main open problem about \LST\ is the following:

\begin{conjecture}\label{cstokes}
The $L^2$--Stokes theorem holds for complex projective varieties.
\end{conjecture}

Conjectures closely related to Conjecture
\ref{cstokes}\ were stated in the basic paper
\cite{CheGorMac:LCIHSAV}, but Conjecture 
\ref{cstokes}\ was not stated explicitly. 
For complex surfaces, Conjecture \ref{cstokes}\ was formulated in \cite{Nag:RLCSAS}.

We now review some results about the validity of \LST.

\begin{enumerate}
\item {\em Complete manifolds:} Gaffney \cite{Gaf:HOEDF},
\cite{Gaf:SSTCRM} showed that $\Delta$ is essentially self-adjoint,
in particular \LST\ holds. 
\item {\em Cones, horns and pseudomanifolds:} Let $(N^n,g_N)$ be
a Riemannian manifold satisfying \LST, and such that $\text{Range } d_{n/2-1}$ is
closed (e.g.\ $N$ compact).
Define the cone ($\gamma=1$)
or horn ($\gamma>1$) over $N$ by
 $$ M=(0,\infty)\times N,\quad g=dx^2 + x^{2\gamma} g_N.$$ 
Cheeger \cite[Thm.\ 2.2]{Che:HTRP} proved that  \LST\ holds for $(M,g)$
if and only if there are no square-integrable
 $n/2$-forms $\alpha$ on $N$ satisfying $d\alpha=d^t\alpha=0$
(if $N$ is compact, this is the cohomological condition  $H^{n/2}(N,\R)=0$).
Using this inductively, he showed \LST\ for 
'admissible Riemannian pseudomanifolds' (Proof of Thm.\ 5.1 in loc.\ cit.).
Youssin \cite{You:LCCH} generalized these results to the $L^p$ Stokes
theorem for any $p\in (1,\infty)$ and carried out a detailed study of their
relation to $L^p$-cohomology.
\item {\em (Conformal) cones, complex projective algebraic curves:}
A different proof for cones over compact $N$ and a generalization to conformal
cones was given by Br\"uning and Lesch \cite{BruLes:KHTCCC}. 
In particular, this implies
\LST\ for  complex projective algebraic curves. The last result was
obtained before by Nagase \cite{Nag:HTSAC}.
\item {\em Functions on real or complex projective algebraic varieties
 with singularities of real codimension at least two:} 
Li and Tian \cite{LiTia:HKBMAV} proved that $d^t_{0,\min}d_{0,\min}$ is 
self-adjoint for such $M$. So, in view of Proposition
\plref{prop2.3} uniqueness holds for $k=0$. This was proved before
by Nagase \cite{Nag:HONSAS}\ and Pati \cite{Pat:HTSAT}\ for complex
surfaces and threefolds, respectively.

\item {\em Complex projective algebraic varieties with isolated 
singularities:} Pardon and Stern \cite{ParSte:PHSLCVIS}\ proved that
the cohomology of the $d_{k,\min}$ and $d_{k,\max}$--complexes
coincide in degrees $k$ with 
$|k-n| \geq 2$, $n=\dim_\C M$.
We will prove below (Propositions \plref{3.1} and \plref{3.2})
that their estimates actually imply the stronger
 $\dmaxdmin{k}$ for $k\not=n-1,n$.

\item {\em Real analytic surfaces with isolated singularities:}
In this case, uniqueness was proved in all degrees by
Grieser \cite{Gri:LGSRAS}.

\item {\em Orbit spaces:} 
If $G$ is a compact Lie group acting isometrically on a smooth compact
Riemannian manifold $X$ then the quotient $X/G$ is a stratified space
with a Riemannian metric on its smooth part $M$ (see \cite{Sja:LCOS}).
Then \LST\ holds on $M$ if and only if $X/G$ is a Witt space, i.e.
the links of all the odd-codimensional strata have vanishing middle
$L^2$-cohomology. A study of when this happens was made for $S^1$-actions
by Sjamaar \cite{Sja:LCOS}.
\end{enumerate} 

The methods used to prove \LST\ all proceed essentially in the following way:
First, one derives certain boundedness or decay properties for forms
 $\omega\in\dom(d_{\max})$, then one uses these to show that, for
a suitable sequence of cutoff-functions
$\{\phi_n\}\in C_0^\infty(M)$, the forms $\phi_n\omega$ converge
to $\omega$ in graph norm. By Lemma \plref{lemma2.1}, 
this implies $\omega\in\dom(d_{\min})$. 
In the first step, it may be useful to
restrict to certain dense (in graph norm) subspaces of $\dom(d_{\max})$ with
better (or easier to derive) properties. 
For example, from the proof of Proposition \plref{prop2.3} (1) 
we see in particular that
$$\dom(d_{k,\max}) \cap \dom(d^t_{k-1,\max}) \cap C^\infty(\Wedge^k M)$$
is a core for $d_{k,\max}$. 
 Thus, we may assume 
\begin{equation} \label{eqomDom}
\omega\in\Ltwo,\quad (d+d^t)\omega\in\Ltwo.
\end{equation}
This allows to use elliptic techniques to
estimate $\omega$.

On functions, one can do slightly better:

\begin{lemma}\label{S-2.5}{\rm (cf. \cite[Sec. 4]{LiTia:HKBMAV}) }Let $(M,g)$
be a Riemannian manifold. Then 
$$\{f\in\cd(d_{0,\max})\, |\,
\sup\limits_{x\in M}|f(x)|<\infty\}$$
is a core for $d_{0,\max}$.
\end{lemma}
\begin{proof} For $f\in\cd(d_{0,\max})\cap\cinf{M}$ choose a sequence $(a_n)$
of regular values of $f$ with $\lim_{n\to\infty} a_n=\infty$
and put $f_n:=\max(-a_n,\min(f(x),a_n))$. Then $f_n$ is bounded
and it is straightforward to check that $f_n\to f$ in the graph
norm of $d$.\end{proof}

In general, if $\omega\in\dom(d_{k,\max})$ and 
$(\phi_n) \subset C_0^\infty(M)$,
in order to prove $\omega\in \dom(d_{k,\min})$ it suffices to show
\begin{align}
&\phi_n\omega  \longrightarrow  \omega \quad\text{ in } L^2(\Wedge^k M),
                                           \label{eqomconv} \\
&\|d(\phi_n\omega) \|_{L^2(\Wedge^k (M))} \leq  C. \label{eqdomconv}
\end{align}
This follows from Lemma \plref{lemma2.1}\ and the analogue of estimate
\eqref{domest}.

By Lebesgue's dominated convergence theorem,
(\plref{eqomconv}) is fulfilled if
$$ \text{Condition 1: } \phi_n \longrightarrow 1 \text{ pointwise and }
           \sup_{x,n} |\phi_n(x)| < \infty.$$
(The second condition could be weakened to $\sup_n |\phi_n(x)| \leq \phi(x)$
for some $\phi$ with $\phi\omega\in \Ltwo$.)
Then, (\plref{eqdomconv}) will be satisfied if 
$$ \text{Condition 2: } \|d\phi_n\wedge\omega\|_{L^2} \leq C$$
(plus $\phi\, d\omega\in\Ltwo$ in case of the weaker Condition 1).
 
Specifically, the proofs of the results cited above proceed in
the following way:
\begin{enumerate}
\item Gaffney finds $\phi_n$ satisfying
Condition 1 and with $d\phi_n$ uniformly bounded as $n\to\infty$.
\item Cheeger proves (\plref{eq6}) directly, using the weak
Hodge decomposition of $\omega\in\dom(d_{\max})$ and estimating its parts
separately. Since validity of \LST\ is a local property (see
\cite[Lemma 4.1]{Che:HTRP}) the result on pseudomanifolds follows by
localization and iteration, using a partition of unity with uniformly
bounded differentials.
\item Br\"uning and Lesch
 show that (\plref{eqomDom}) implies $\omega\in\dom(d_{\min})$ modulo
the pull-back of a harmonic $n/2$-form on $N$, by writing
the equation $D\omega=f\in\Ltwo$ as an operator-valued singular ordinary
differential equation in the axis variable $x$, whose solutions can be 
analyzed fairly explicitly.
\item
Li and Tian first restrict to bounded functions using Lemma \ref{S-2.5}
and then use a sequence $\phi_n$ satisfying
Condition 1 and $\|d\phi_n\|_{\Ltwo} \to 0$, which clearly implies
Condition 2. 
\item
Pardon and Stern show, using the K\"ahler identities,
 that $d\omega,d^t\omega\in\Ltwo$ implies $\omega/r\in\Ltwo$ if $|k-n|\geq 2$,
where $r$ is the distance from the singularity 
\cite[Prop.\ 2.27]{ParSte:PHSLCVIS}. 
Then they use a 
cutoff sequence satisfying Condition 1 and with $rd\phi_n$ tending to zero
uniformly as $n\to\infty$. This again implies Condition 2.
\item Grieser shows that a neighborhood of a singular point is 
quasi-isometric to a union of cones and horns, so the result follows
from (2).
\item This follows from the results on Riemannian pseudomanifolds of (2)
since the metric on $M$ is locally quasi-isometric to a metric of
iterated cone type.

\end{enumerate}

We now show that on real algebraic varieties \LST\ may very well
fail, even near isolated singularities. For this, we construct
examples which locally near the singularity are quasi-isometric to
\begin{equation}\begin{split}
M & =  \R_+ \times N_1 \times N_2 \\
g & =  dr^2 + r^{2\alpha_1} g_1 + r^{2\alpha_2} g_2 
                \end{split}
  \label{G2.8}
\end{equation}
with compact oriented Riemannian manifolds $(N_i,g_i)$ and
real numbers $\alpha_i\geq 1$, $i=1,2$. Here, $r$ denotes the first variable.

\begin{lemma}\label{llocalcalc}
Let $(M,g)$ be as in \myref{G2.8}. 
If $H^k(N_2)\not= 0$ then \LST\ fails for $k$-forms
supported near $r=0$, where
\begin{equation} \label{eqk}
k = \frac{\alpha_1n_1 + \alpha_2n_2}{2\alpha_2},\quad  n_i=\dim N_i.
\end{equation}
\end{lemma}

\begin{proof}
Let $\omega$ be a non-zero harmonic $k$-form on $N_2$, i.e.\
$d_2\omega = d^t_2\omega = 0$ (in this proof, the subscripts 1,2 refer
to $N_1,N_2$ respectively). Choose cut-off functions $\phi,\psi\in
C_0^\infty([0,\infty))$ such that $\phi=1$ near zero, $\supp \psi \subset
\{\phi=1\}$, and $\psi(0)\not=0$.
For simplicity, denote the pullbacks of $\omega,\phi,\psi$ to $M$ by the same
letters. 

Set $\mu = \phi\omega\in\Omega^k(M)$, $ \nu = \psi
dr\wedge\omega\in\Omega^{k+1}(M)$. 
We claim that $$ (d\mu,\nu) \not= (\mu,d^t\nu).$$
To prove this, recall $(d\mu,\nu) - (\mu,d^t\nu) = \int_M d(m\wedge *\nu)$.
A simple calculation shows
$*\nu = \pm \psi(r) r^\beta dvol_1 \wedge *_2\omega$ (the sign is unimportant),
 where
$\beta = \alpha_1n_1 + \alpha_2(n_2-2k)$.
So $\mu\wedge *\nu = \pm \psi(r) r^\beta dvol_1 \wedge (\omega\wedge
*_2\omega)$. Since $\omega$ is closed and coclosed, we get
$d(\mu\wedge *\nu) = \pm d(\psi r^\beta) \wedge dvol_1 \wedge(\omega\wedge
*_2\omega)$.
Therefore, if $\beta\geq0$ we have
$\int_M d(\mu\wedge *\nu) = \pm \|\omega\|^2_{L^2(N_2)} \psi(r) r^\beta_{|r=0}$.
This is non-zero for $\beta=0$, i.e.\ $k$ as in (\plref{eqk}).
\end{proof}

For example, if $$ \alpha_1n_1 = \alpha_2 n_2$$
then $k=n_2$, so \LST\ fails since always $H^{n_2}(N_2)\not=0$. By symmetry,
\LST\ also fails for $k=n_1$.

\begin{lemma}\label{lrealex}
Let $n,m,p,q$ be positive integers, $p>q$. For $x\in\R^n, y\in\R^m,z\in\R$
define the homogeneous polynomial $f$ by
$$f(x,y,z) = |x|^{2p} - |y|^{2q} z^{2(p-q)}$$
and the variety $V$ by 
$$ V = \{[x,y,z] \,|\, f(x,y,z) = 0\} \subset \RP^{n+m}.$$
Then, the singularities of $V$ are given by 
\begin{equation}
\sing V=\{[0,0,1]\}\cup \{[0,y,0]\,|\, y\in\R^m\setminus\{0\}\}.
\label{G2.9}
\end{equation}

Thus, $V$ has an isolated singularity at $P=[0,0,1]$ and one
singular stratum isomorphic to $\RP^{m-1}$. 
Given any Riemannian
metric on $\RP^{n+m}$, a (pointed)
neighborhood of $P$ in $V$ is quasi-isometric to
a neighborhood of $r=0$ in 
\begin{equation}\begin{split}
M&= \R_+ \times S^{n-1} \times S^{m-1}  \\
g& =  dr^2 + r^2 g_{S^{n-1}} + r^{2\alpha} g_{S^{m-1}}
                \end{split}
             \label{eqmetric}
\end{equation}
where $\alpha = p/q$.

\end{lemma}

\begin{proof} That $\sing V$ is given by the right hand side of
\myref{G2.9} can be checked by a straightforward calculation.

For $z\not=0$, we use the standard chart $(x,y)\in\R^{n+m} \to
[x,y,1]\in\RP^{n+m}$. In this chart, $V$ is given by $\{|x|^{2p} = |y|^{2q}\}$,
which clearly has an isolated singularity at $(0,0)$.
Parametrize $V$ near but outside $(0,0)$ by
\begin{align*}
\phi: &\R_+\times S^{n-1} \times S^{m-1} \longrightarrow V  \\
&(r,v,w)  \longmapsto  (rv,r^\alpha w),\quad \alpha = p/q.
\end{align*}
Any Riemannian metric on $\RP^{n+m}$ is quasi-isometric to the standard
Euclidean metric $|dx|^2 + |dy|^2$ near $P$. By the usual formula for
polar coordinates, we have
\begin{align*}
\phi^*(|dx|^2)& = dr^2 + r^2 g_{S^{n-1}} \\
\phi^*(|dy|^2)& = (d(r^\alpha))^2 + r^{2\alpha} g_{S^{m-1}},
\end{align*}
so $\phi^*(|dx|^2+|dy|^2) = (1+\alpha^2r^{2(\alpha-1)}) dr^2 + r^2 g_{S^{n-1}}
+ r^{2\alpha} g_{S^{m-1}}$, which is quasi-isometric to (\plref{eqmetric})
for bounded $r$. \end{proof}

\begin{prop}\label{prealex}
Let $V$ be as above and assume $p(m-1)=q(n-1)$.
Then \LST\ fails for $V$ for $(n-1)$-forms and for
$(m-1)$-forms.
\end{prop}

\begin{proof}
This follows from Lemma \plref{lrealex}\ and
Lemma \plref{llocalcalc}\ in the special case
mentioned after its proof.
\end{proof}

In the previous proposition we may choose
$m$ odd and $n$ even (e.g. $m=2k+1, n=2k+2, p=2k+1, q=2k$, $k>0$).
Then $V$ is even--dimensional and all strata are of even codimension.

We can even slightly modifiy $V$ to produce a variety
with one isolated singularity for which \LST\ does
not hold:

\begin{lemma}\label{2.9} Let $f,g\in\cinf{\R}$ with $f^{(j)}(0)=0,$
$0\le j< k,$ $f^{(k)}(0)\not=0$, $g^{(j)}(0)=0,$
$0\le j< l,$ $g^{(l)}(0)\not=0$, $k\ge l$. Then, locally
near $0$, the quasi--isometry class of 
$$V_{f,g}:=\{(x,y)\in\R^{n+m}\,|\, f(|x|^2)=g(|y|^2)\}$$
depends only on $k$ and $l$.
\end{lemma}
\begin{proof} 
W.l.o.g. we may assume that $f$ and $g$ are non-negative for small positive
argument. Let $\tilde f(r)=r^k, \tilde g(r)=r^l$. 
Then, there exist $\varphi,\psi\in\cinf{\R}$, 
$\varphi(0)=\psi(0)=0,$  $\varphi'(t),\psi'(t)>0$ for $|t|<\eps$,
such that
$$\tilde f=f\circ\varphi,\quad \tilde g=g\circ\psi.$$
We put $\varphi_1(t):=\varphi(t)/t, \psi_1(t):=\psi(t)/t$ and
$$\Phi(x,y):= (\varphi_1(|x|^2)^{1/2} x, \psi_1(|y|^2)^{1/2} y),
   \quad (x,y)\in\R^{n+m},\quad |x|,|y|<\eps.$$
Then $\Phi$ is a local diffeomorphism at $(0,0)$ mapping
$V_{\tilde f,\tilde g}\cap U$ onto $V_{f,g}\cap U$ for
some neighborhood $U$ of $(0,0)$. Hence
$V_{\tilde f,\tilde g}\cap U$ and $V_{f,g}\cap U$ are quasi--isometric.
\end{proof}

Summing up, we have:

\begin{prop}\label{2.10} Under the assumptions of Lemma \plref{lrealex},
consider the variety
$$ W = \{[x,y,z] \,|\, |x|^{2p}=|y|^{2q} z^{2(p-q)}+|y|^{2p}\}
   \subset \RP^{n+m}.$$
Then $\sing W=\{[0,0,1]\},$ i.e. $W$ has exactly one singular point.
Moreover, given any Riemannian
metric on $\RP^{n+m}$, a (pointed)
neighborhood of $[0,0,1]$ in $W$ is quasi-isometric to
a pointed neighborhood of $[0,0,1]$ in the variety $V$ of
Lemma \plref{lrealex}. 

In particular, if  $p(m-1)=q(n-1)$, then \LST\ fails to hold
for $W$.
\end{prop}
\begin{proof} The determination of $\sing W$ is straightforward. The other
statements follow from the previous discussion.\end{proof}


\section{Complex algebraic varieties}
In this section we discuss the $L^2$--Stokes theorem and the
discreteness of the spectrum of the Laplace-Beltrami operator
on complex projective varieties with isolated
singularities. The main ingredients
are an estimate due to Pardon and Stern \cite{ParSte:PHSLCVIS}\
and results by Br\"uning and Lesch \cite{BruLes:KHTCCC}.

Throughout this section, $V\subset \CP^N$ denotes an algebraic variety
with isolated singularities and $M=V\setminus \sing V$ its smooth 
part, $n:= \dim_{\C} M$. Let $g$ be a K\"ahler metric on $M$ which is
induced by a K\"ahler metric on $\CP^N$.
Fix a smooth function $r : M \to (0,\infty)$ which, near any singularity
$p\in \sing V$, is comparable to the distance from $p$.

Denote by $\pl_{p,q},  \plbar_{p,q}$ the Dolbeault operators on
forms of type $(p,q)$.
Pardon and Stern show in \cite{ParSte:PHSLCVIS}:

\begin{prop}\label{3.1}
\begin{enumerate}
\item Assume $n-p-q \geq 2$. If $\omega\in \dom(\pl_{p,q,\max})
     \cap \dom (\pl_{p-1,q,\max}^t)$ then
     $$\omega /r \in L^2(\Wedge^{p,q}M)$$
     and $\omega \in \dom((\plbar + \plbar^t)_{\min}) 
     \cap \dom((\pl+\pl^t)_{\min})$.
\item Assume $n-k\geq2$. If $\omega \in \dom (d_{k,\max}) \cap
     \dom(d_{k-1,\max}^t)$ then
     $$ \omega /r \in L^2(\Wedge^kM)$$
     and $\omega \in \dom ((d+d^t)_{\min})$.
\end{enumerate}
\end{prop}

For the proof see \cite{ParSte:PHSLCVIS}, Prop. 2.27 and Lemma 2.18.

We combine this with a result in \cite{BruLes:KHTCCC}.

\begin{prop}\label{3.2}
Let $(C_0^\infty(E),d)$ be an elliptic complex. If for some $k\geq 0$
\begin{equation} \label{eqDmaxmin}
 \dom (d_{k,\max}) \cap \dom (d^t_{k-1,\max}) 
    \subset \dom((d+d^t)_{\min}) 
\end{equation}    
then we have
\begin{align}
d_{k,\max} & = d_{k,\min}, \label{4a} \\
d^t_{k-1,\max} & = d^t_{k-1,\min}, \label{4b}
\end{align}
and
$$ \Delta^{r/a}_k = d_{k-1,\min} d^t_{k-1,\min} + 
                    d^t_{k,\min} d_{k,\min}  = \Delta^{\cF}_k,$$
where $\Delta^{\cF}_k$ denotes the Friedrichs extension of the 
Laplacian $\Delta_k$.
\end{prop}

For the proof see \cite{BruLes:KHTCCC}, Lemma 3.3.
Note that (\ref{eqDmaxmin}) is equivalent to
$$ D_{k,\max} = D_{k,\min}$$
for $D_k := d_k+d^t_{k-1}$. However, one needs to be careful when considering
the 'rolled-up' complex, i.e.\ the operator $D=\bigoplus_{k=0}^N D_k$:
Even if (\ref{eqDmaxmin}) holds for all $k$, this does \emph{not} imply
$D_{\max}=D_{\min}$! The reason for this is that, in general,
$$ \dom(D_{\max}) \supsetneqq \bigoplus_k \dom(D_{\max}) \cap L^2(\Wedge^k M)
   $$
since the ranges of $d_{\max}$ and $d^t_{\max}$ might not be orthogonal.   
An example can be easily constructed using the de Rham complex on the unit
disk in $\R^2$.
\vspace{2mm}

\noindent
{\em Proof of Theorem \ref{3.3}.}  
First, consider the de Rham complex.
Propositions \plref{3.1}\ and \plref{3.2}\ imply that \myref{4a}\ and
\myref{4b}\ hold for $k=0,\ldots,n-2$. 
From $d_k = \pm * d^t_{2n-k-1} *$ we get the same relations for
$k=n+2,\ldots,2n$. Finally, \myref{4b}, applied with $k=n+2$,
gives \myref{4a}\ for $k=n+1$ by taking adjoints.
This proves \myref{3a}. The same argument applied to the
Dolbeault complex proves \myref{3b}.

The first and third equality in \myref{4}\ follow from the last identity
in Proposition \plref{3.2}, applied to the de Rham and Dolbeault complex,
respectively.
It remains to prove 
\begin{equation}\label{friedrich}
 \Delta^{\cF}_k  = 2 \bigoplus_{p+q=k} \Delta^{\cF}_{p,q,\pl}.
\end{equation}
Now the K\"ahler identities (see \cite{Wel:DACM}, Ch.\ V, Thm.\ 4.7)
show that, for a compactly supported smooth $k$-form
$\phi = \bigoplus_{p+q=k} \phi_{p,q}$, one has
\begin{equation}\label{dirsum}
\Delta_k\phi = 2\bigoplus_{p+q=k} \Delta_{p,q,\partial}\, \phi_{p,q}.
\end{equation}
The domain of the Friedrichs extension of $\Delta_k$
is defined as the completion of
$C_0^\infty(\Wedge^k M)$ with respect to the scalar product
$$ (\phi,\psi)_{\Delta_k^{\cF}} = (\Delta_k \phi,\psi) + (\phi,\psi).$$
By \eqref{dirsum} this scalar product is the direct sum of the
analogous scalar products $(\ ,\ )_{\Delta_{p,q}^{\cF}}$ on
$C_0^\infty(\Wedge^{p,q}M)$,  $p+q=k$, 
so we obtain \eqref{friedrich}.
\qed

\remark For the operators $\pl,\plbar$ this result is sharp. This can be
seen already in the case of algebraic curves: In \cite{BruPeySchro:OAC}\
it was shown that for $n=1$ in general
$\pl_{0,\max} \not= \pl_{0,\min}, \pl_{0,1,\max} \not= \pl_{0,1,\min}$.

We now turn to the question of discreteness. Recall that a
self--adjoint operator
$T$ in some Hilbert space $H$ is called {\it discrete} if
its spectrum consists only of eigenvalues of finite multiplicity
(with $\infty$ as the only accumulation point). Note that
$T$ is discrete if and only if it has a compact resolvent
or equivalently if the embedding $\cd(T)\hookrightarrow H$ is compact.
Here, $\cd(T)$ carries the graph topology.
\begin{prop}\label{3.4}
Let $(M,g)$ be a Riemannian manifold such that there is a 
function $r\in\cinf{M}$, $r>0$, with the following properties:
\begin{enumerate}\labelroman
\item For $\eps>0$ the set $r^{-1}([\eps,\infty))$ is compact,
\item $f\in\dom(\Delta_k^\cF)$ implies $f/r \in L^2(\Wedge^k (M))$.
\end{enumerate}
Then $\Delta_k^\cF$ is discrete.

In particular, if $M$ is as in Theorem \plref{3.3}\ then $\Delta_k^\cF$
is discrete for $k\not= n,n\pm 1$.
\end{prop}

\begin{proof}  \newcommand{\dfk}{\Delta_k^\cF}
We first note that (ii) implies
$$\|f/r\| \leq C \|(I+\dfk)f\| $$
for $f\in\dom(\dfk)$. ($\|\cdot\|$ denotes the $L^2$ norm.)
This follows from the closed graph theorem since the operator of
multiplication by $1/r$ is closable with domain containing $\dom(\dfk)$.
In order to prove discreteness, it is enough to show that the embedding
$\dom(\dfk) \hookrightarrow L^2$ is compact, where $\dom(\dfk)$ carries the graph norm.
Thus, let $(\phi_n) \subset \dom(\dfk)$ be a bounded sequence, i.e.
$$ \|\phi_n\|, \| \dfk \phi_n\| \leq C_1.$$
We need to show that $(\phi_n)$ has a  subsequence that is convergent
in $L^2(\Wedge^k(M))$.

Since $\Delta_k$ is an elliptic differential operator and since $\{r\geq1\}$
is compact there exists a subsequence $(\phi_n^{(1)})$ which converges in
$L^2(\Wedge^k(\{r\geq1\}))$. 
Choose a subsequence $(\phi_n^{(2)})$ of $(\phi_n^{(1)})$
which converges in $L^2(\Wedge^k(\{r\geq 1/2\}))$. 
Continuing in this way and then
using a diagonal argument, we find a subsequence $(\phi_n^{(\infty)})$ such
that
$$ (\phi^{(\infty)}_n |\{r\geq 1/m\}) \text{ converges in }
    L^2(\Wedge^k(\{r\geq 1/m\})), \quad \text{ for all } m. $$
We show that $(\phi_n^{(\infty)})$ is a Cauchy sequence in $L^2(\Wedge^k(M))$:
First, we estimate
$$ \|\phi_n^{(\infty)}\|_{L^2(\{r\leq 1/m\})}
    \leq \frac1m \| \frac1r \phi_n^{(\infty)}\|
    \leq \frac{C}m \|(I+\dfk)\phi_n^{(\infty)}\| \leq C'/m.  $$
Given $\eps>0$ choose $m$ such that $C'/m < \eps/3$. Then
$$ \|\phi_k^{(\infty)} - \phi_l^{(\infty)} \|_{L^2(\Wedge^k(M))}
\leq \frac23 \eps + \|\phi_k^{(\infty)} - \phi_l^{(\infty)} \|_{L^2(\{r\geq 1/m\})}.$$
Since the last term tends to zero as $k,l\to\infty$, we are done.
\end{proof}

\noindent {\em Proof of Theorem \ref{3.5}.}
This follows by applying Proposition \plref{3.4} to the situation of Theorem
\plref{3.3} in view of Proposition \plref{3.1}.
\qed

\addcontentsline{toc}{section}{References}
\renewcommand\MakeUppercase{\relax}
\renewcommand\refname{{\upshape References}}

\begin{thebibliography}{10}

\bibitem{BruLes:HC}
\textsc{J.~Br{\"u}ning {\upshape and} M.~Lesch}: \emph{Hilbert complexes}.
\newblock J. Funct. Anal. \textbf{108} (1992), 88--132

\bibitem{BruLes:KHTCCC}
\textsc{J.~Br{\"u}ning {\upshape and} M.~Lesch}: \emph{K\"ahler--Hodge theory
  for conformal complex cones}.
\newblock Geom. Funct. Anal. \textbf{3} (1993), 439--473

\bibitem{BruLes:SGAC}
\textsc{J.~Br\"uning {\upshape and} M.~Lesch}: \emph{On the spectral geometry
  of algebraic curves}.
\newblock J. reine angew. Math. \textbf{474} (1996), 25--66

\bibitem{BruPeySchro:OAC}
\textsc{J.~Br\"uning, N.~Peyerimhoff, {\upshape and} H.~Schr\"oder}: \emph{The
  $\bar\partial$--operator on algebraic curves}.
\newblock Commun. Math. Phys. \textbf{129} (1990), 525--534

\bibitem{Che:HTRP}
\textsc{J.~Cheeger}: \emph{On the Hodge theory of Riemannian pseudomanifolds}.
\newblock Proc. Symp. Pure Math. \textbf{36} (1980), 91--146

\bibitem{Che:HTCC}
\textsc{J.~Cheeger}: \emph{Hodge theory of complex cones}.
\newblock Ast\'erisque \textbf{102} (1983), 118--134

\bibitem{Che:SGSRS}
\textsc{J.~Cheeger}: \emph{Spectral geometry of singular Riemannian spaces}.
\newblock J. Diff. Geom. \textbf{18} (1983), 575--657

\bibitem{CheGorMac:LCIHSAV}
\textsc{J.~Cheeger, M.~Goresky, {\upshape and} R.~MacPherson}:
  \emph{$L^2$--cohomology and intersection homology of singular algebraic
  varieties}.
\newblock Annals of Math. Studies \textbf{102} (1982), 303--340.
\newblock Princeton University Press, Princeton, NJ

\bibitem{Gaf:HOEDF}
\textsc{M.~P. Gaffney}: \emph{The harmonic operator for exterior differential
  forms}.
\newblock Proc. Nat. Acad. Sci. USA \textbf{37} (1951), 48--50

\bibitem{Gaf:SSTCRM}
\textsc{M.~P. Gaffney}: \emph{A special Stokes theorem for complete Riemannian
  manifolds}.
\newblock Ann. of Math., II. Ser. \textbf{60} (1954), 140--145

\bibitem{Gri:LGSRAS}
\textsc{D.~Grieser}: \emph{Local geometry of singular real analytic surfaces}.
\newblock Preprint Erwin-Schr\"odinger-Institut Wien

\bibitem{HsiPat:LCNAS}
\textsc{W.~C. Hsiang {\upshape and} V.~Pati}: \emph{$L\sp 2$-cohomology of
  normal algebraic surfaces. I}.
\newblock Invent. Math. \textbf{81} (1985), 395--412

\bibitem{LiTia:HKBMAV}
\textsc{P.~Li {\upshape and} G.~Tian}: \emph{On the heat kernel of the Bergmann
  metric on algebraic varieties}.
\newblock J. Amer. Math. Soc. \textbf{8} (1995), 857--877

\bibitem{Mel:APSIT}
\textsc{R.~B. Melrose}: \emph{The Atiyah--Patodi--Singer index theorem}.
\newblock A K Peters, Wellesley, Massachusetts (1993)

\bibitem{Nag:HONSAS}
\textsc{M.~Nagase}: \emph{On the heat operators of normal singular algebraic
  surfaces}.
\newblock J. Diff. Geom. \textbf{28} (1988), 37--57

\bibitem{Nag:RLCSAS}
\textsc{M.~Nagase}: \emph{Remarks on the $L\sp 2$-cohomology of singular
  algebraic surfaces}.
\newblock J. Math. Soc. Japan \textbf{41} (1989), 97--116

\bibitem{Nag:HTSAC}
\textsc{M.~Nagase}: \emph{Hodge theory of singular algebraic curves}.
\newblock Proc. Amer. Math. Soc. \textbf{108} (1990), 1095--1101

\bibitem{Ohs:LCCS}
\textsc{T.~Ohsawa}: \emph{On the $L\sp 2$ cohomology of complex spaces}.
\newblock Math. Z. \textbf{209} (1992), 519--530

\bibitem{Ohs:LCCSII}
\textsc{T.~Ohsawa}: \emph{On the $L\sp 2$ cohomology of complex spaces. II}.
\newblock Nagoya Math. J. \textbf{127} (1992), 39--59

\bibitem{ParSte:PHSLCVIS}
\textsc{W.~Pardon {\upshape and} M.~Stern}: \emph{Pure Hodge structure on the
  $L_2$-cohomology of varieties with isolated singularities}.
\newblock J.\ reine angew.\ Math.\ \textbf{533} (2001), 55-80

\bibitem{ParSte:LCCPV}
\textsc{W.~Pardon {\upshape and} M.~Stern}: \emph{$L\sp
  2-\bar\partial$-cohomology of complex projective varieties}.
\newblock J. Am. Math. Soc. \textbf{4} (1991), 603--621

\bibitem{Pat:HTSAT}
\textsc{V.~Pati}: \emph{The heat trace on singular algebraic threefolds}.
\newblock J. Differ. Geom. \textbf{37} (1993), 245--261

\bibitem{Pat:LATIS}
\textsc{V.~Pati}: \emph{The Laplacian on algebraic threefolds with isolated
  singularities}.
\newblock Proc. Indian Acad. Sci., Math. Sci. \textbf{104} (1994), 435--481

\bibitem{RemSchul:ITEBP}
\textsc{S. Rempel {\upshape and} B.-W. Schulze}: 
\emph{Index Theory of Elliptic Boundary Value Problems}.
\newblock Akademie Verlag, Berlin (1982)

\bibitem{Schul:ECMWCS}
\textsc{B.-W. Schulze}: \emph{Elliptic Complexes on Manifolds with Conical 
Singularities}.
\newblock Teubner Texte zur Math. (Leipzig)
\textbf{106} (1988), pp. 170--223.


\bibitem{Schul:PDBVPCSA}
\textsc{B.-W. Schulze}: \emph{Pseudo-Differential Boundary Value Problems,
  Conical Singularities and Asymptotics}.
\newblock Akademie Verlag, Berlin (1994)

\bibitem{Seg:FC}
\textsc{G. Segal}: \emph{Fredholm complexes}.
\newblock Quarterly J. Math. \textbf{21} (1970), 385--402.

\bibitem{Sja:LCOS}
\textsc{R.~Sjamaar}: \emph{{$L^2$}-cohomology of orbit spaces}.
\newblock Topology Appl. \textbf{45} (1992), 1--11

\bibitem{Tay:PO}
\textsc{M.~Taylor}: \emph{Pseudodifferential Operators}.
\newblock Princeton University Press (1981)

\bibitem{Wel:DACM}
\textsc{R.~O.~Wells}: \emph{Differential Analysis on Complex Manifolds}.
\newblock (2nd ed.) Springer Verlag (1980), New York.

\bibitem{You:LCCH}
\textsc{B.~Youssin}: \emph{{$L^p$} cohomology of cones and horns}.
\newblock J. Differential Geom. \textbf{39} (1994), 559--603

\end{thebibliography}

\signature
\end{document}